\documentclass[a4paper]{article}
\usepackage[utf8]{inputenc}
\usepackage{amsmath,amsthm,amssymb,amsfonts}
\usepackage[german,english]{babel}

\theoremstyle{plain}
\newtheorem{theorem}{Theorem}[section]
\newtheorem{lemma}[theorem]{Lemma}
\newtheorem{claim}[theorem]{Claim}
\newtheorem{corol}[theorem]{Corollary}
\newtheorem*{remark}{Remark}

\title {
    Maximum number of symmetric extensions in the random graph
}
\author {\Large
    S. Vakhrushev\footnote{Department of Mathematics and Computer Science, Saint-Petersburg University, Saint Petersburg, Russia}, M. Zhukovskii\footnote{Department of Computer Science, The University of Sheffield, Sheffield S1 4DP, United Kingdom}
}
\date{}

\usepackage[]{hyperref}

\def\BN{\mathbb{N}}

\def\BR{\mathbb{R}}

\usepackage{subfigure} 

\newcounter{image}
\newenvironment{image}[1][]{\refstepcounter{image} 
	\small Fig.~\theimage. #1 \rmfamily}{}

\usepackage{pgfplots}

\newenvironment{Proof1} 
{\par\noindent{\it Proof for the first case:}} 
{\hfill$\scriptstyle\blacksquare$} 

\newenvironment{Proof2} 
{\par\noindent{\it Proof for the second case:}} 
{\hfill$\scriptstyle\blacksquare$} 
\begin{document}

\maketitle

\renewcommand{\abstractname}{Abstract}
\begin{abstract}
It is known that after an appropriate rescaling the maximum degree of the binomial random graph converges in distribution to a Gumbel random variable. The same holds true for the maximum number of common neighbours of a $k$-vertex set, and for the maximum number of $s$-cliques sharing a single vertex. Can these results be generalised to the maximum number of extensions of a $k$-vertex set for any given way of extending of a $k$-vertex set by an $s$-vertex set? In this paper, we generalise the above mentioned results to a class of ``symmetric extensions'' and show that the limit distribution is not necessarily from the Gumbel family.
\end{abstract}

\section{Introduction}
Bollob\'as~\cite{bolobas} and Ivchenko~\cite{ivchenko}  proved that under some restrictions on the edge probability $p$, the (appropriately rescaled) maximum degree converges in distribution to a Gumbel random variable.

\begin{theorem}[B. Bollob\'as~\cite{bolobas}]
Let $p=\mathrm{const}\in(0,1)$. Let $\Delta_n$ be the maximum degree of $G(n, p)$. For every integer $n \geq 2$ set
$$
a_n = pn + \sqrt{2p(1-p)n \ln n}\left(1 - \frac{\ln \ln n}{4\ln n} - \frac{\ln (2\sqrt{\pi})}{2 \ln n}\right),\quad 
b_n = \sqrt{\frac{p(1-p)n}{2 \ln n}}.
$$
Then
$$
\frac{\Delta_n - a_n}{b_n}\stackrel{d}\to\eta,\quad n\to\infty,
$$
where $\eta$ has cdf
$e^{-e^{-x}}$ (i.e. it is a standard Gumbel random variable), and $\stackrel{d}\to$ denotes convergence in distribution.
\label{th:max_deg}
\end{theorem}

This result was extended by Ivchenko~\cite{ivchenko} to $p=o(1)$ such that $\frac{pn}{\ln^3 n}\to\infty$. \\

The central result of the extreme value theory is the Fisher--Tippet--Gnedenko theorem \cite{FT,ftg} claiming that, if, for an infinite sequence of independent and identically distributed (i.i.d.) random variables $\{\xi_i\}_{i\in\mathbb{N}}$ and some non-random $a_n, b_n$ the distribution of $\frac{\xi^{(n)} - a_n}{b_n}$ converges weakly to a non-degenerate distribution (here, as usual, $\xi^{(n)}=\max\{\xi_1,\ldots\xi_n\}$), then this limit distribution belongs to one of the following three families of distributions: Gumbel, Weibull or Fr\'{e}chet, and the conditions for the limit distribution to belong to one of theses families are known. Note that this result is not applicable to the degree sequences of random graphs since they constitute triangular arrays of dependent random variables. However, the degree sequence can be approximated by independent binomial random variables in the following sense. A fixed vertex of $G(n, p)$ has the binomial distribution $\mathrm{Bin}(n-1, p)$ with $n-1$ trails and success probability $p$. In \cite{nad_mitov} it was proven that, for the maximum $D_N$ of $N$ independent binomial random variables $\xi_{N,1}, \xi_{N,2} , \ldots , \xi_{N,N}\sim\mathrm{Bin}(M, p)$, where $M = M(N) = \omega(\ln^3 N)$, $p = \mathrm{const}$, and for every $x \in \mathbb{R}$, the following is true: 
\begin{equation}
\nonumber
{\sf Pr}\left(D_N \leq pM + \sqrt{2p(1-p)M \ln N}\left[1 - \frac{\ln \ln N}{4\ln N} - \frac{2\sqrt{\pi}}{2\ln N} + \frac{x}{2\ln N}\right]\right) \rightarrow e^{-e^{-x}} \text{as } N \rightarrow \infty.
\end{equation}
It is easy to see that in the case $M = n - 1, N = n$ this result gives the same scaling constants and limit distribution as in Theorem~\ref{th:max_deg}. This is not unexpected since every pair of vertices in $G(n, p)$ is almost independent --- the dependency is only due to the single adjacency relation between these two vertices. 

However, as we will see below, the limit distributions of similar statistics in $G(n,p)$ not necessarily belong to any of the above three families of distributions.\\

To work with dependent random variables (degrees), Bollob\'as used the method of moments. Namely, let us denote by $X$ the number of vertices with degree greater than $a_n + b_n x$. It turns out that the $r$-th moment of the random variable $X$ converges in distribution to the $r$-th moment of the Poisson random variable with mean $e^{-x}$. From this it follows (see \cite[Theorems 30.1, 30.2]{moments}) that $\lim_{n\to\infty}{\sf Pr} (X = 0) =e^{-e^{-x}}$, which implies the result. \\

Recently~\cite{common}, Rodionov and the second author of the paper generalised Theorem~\ref{th:max_deg} for the maximum number of common neighbours of $k$ vertices $\Delta_{n, k}$ in $G(n, p)$, where $k$ is an arbitrary fixed positive integer. Let $p^k \gg \frac{\ln^3 n}{n},  1 - p \gg \sqrt{\frac{\ln \ln n}{n}}$, then appropriately scaling $\Delta_{n, k}$ converges in distribution to a standard Gumbel random variable as well. The authors used a different approach for the following reasons: 

(1) in the case $k>1$ the variance of the analogous random variables approaches infinity that makes the method of moments no longer applicable directly;

(2) it is computationally difficult (and not clear that it is possible to do in general) to estimate higher moments of the analogous random variable $X$. 

But it turns out that it is enough to condition the probability space on certain ``frequent'' events, then, for the conditional probability, prove that ${\sf E } X (X - 1) \sim ({\sf E }X)^2$, and finally apply some bounds on the probability of ``non-existence'' that are inspired by the method of Arratia et al \cite{arratia}. Note that another possible approach to overcome dependencies between weakly dependent random variables is the Stein--Chen method (see, for example,~\cite{stein-chen}) for establishing Poisson approximations.  For example, Malinovsky~\cite{malinovsky} recently presented a proof of Theorem~\ref{th:max_deg} using this method.

\vspace{\baselineskip}

Finally, in~\cite{main} a similar result for the maximum number of $s$-cliques sharing a single vertex was proven.

\vspace{\baselineskip}

Note that all the above statistics are particular cases of \textit{extension numbers} that were studied by Spencer in~\cite{spencer2, spencer}, who was inspired by the fact that properties of these statistics constitute the basis of the argument for the validity of first order $0$-$1$ laws for sparse random graphs \cite{zero-one-when, zero-one-law}. These statistics also appear to be useful in many other applications, see, e.g.~\cite{rand-proc-trian, rand-proc-dynconc, rand-proc-evol}. An extension is simply a rooted subgraph of a given graph isomorphic to a fixed pattern rooted graph. Formally, let $H$ be a graph with a distinguished \textit{set of roots} $R=\{u_1,\ldots,u_t\}$, and let $S=\{u_{t+1},\ldots,u_s\}$ be all the other vertices of $H$ (\textit{expansion} set). An {\it $(R,H)$-extension} of a tuple of vertices $T=(x_1,\ldots,x_t)$ is a graph $G$ on $\{x_1,\ldots,x_s\}$ such that for all $i<j$ such that $j>t$, the vertices $u_i, u_j$ are adjacent in $H$ if and only if $x_i,x_j$ are adjacent in $G$.  Fix a rooted graph $(R,H)$ and a $t$-tuple $T$ from $[n]:=\{1,\ldots,n\}$. Denote by $X(T):=X_{(R,H)}(T)$ --- \textit{extension count} --- the total number of $(R,H)$-extensions of $T$ in $G(n,p)$ (note that we count extensions as not necessarily induced subgraphs). Spencer~\cite{spencer} proved the law of large numbers for the number of extensions in the case when $(R,H)$ is \textit{grounded} (there is at least one edge between the set of roots and the expansion set in $H$) and \textit{strictly balanced} (extensions in which all proper subextensions have a strictly lower density) rooted graph and $p$ is large enough. 
These results were recently refined in~\cite{warnke}.
\vspace{\baselineskip}

\begin{figure}[htbp]	
	\subfigure[edge extension] 
	{
		\begin{minipage}{2.5cm}
			\centering          
			\includegraphics[scale=0.4]{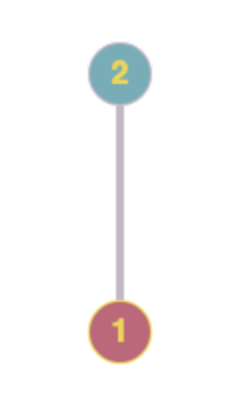}   
		\end{minipage}
	}
	\subfigure[common neighbour extension for $k=5$] 
	{
		\begin{minipage}{5.3cm}
			\centering      
			\includegraphics[scale=0.35]{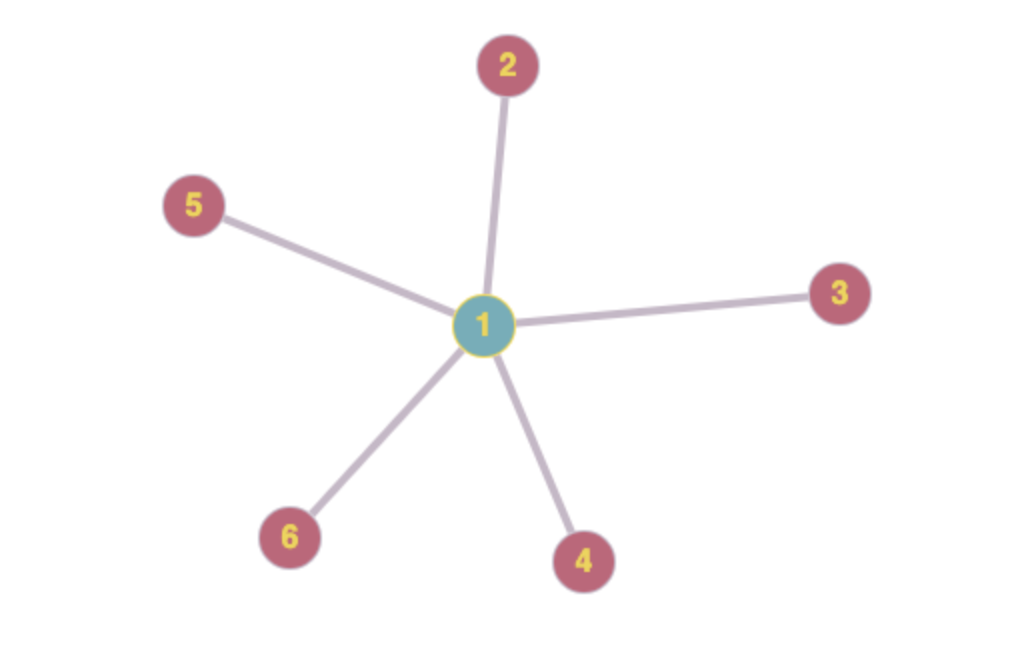}   
		\end{minipage}
	}
	\subfigure[extension of a single vertex by a 5-clique] 
	{
		\begin{minipage}{5.3cm}
			\centering      
			\includegraphics[scale=0.4]{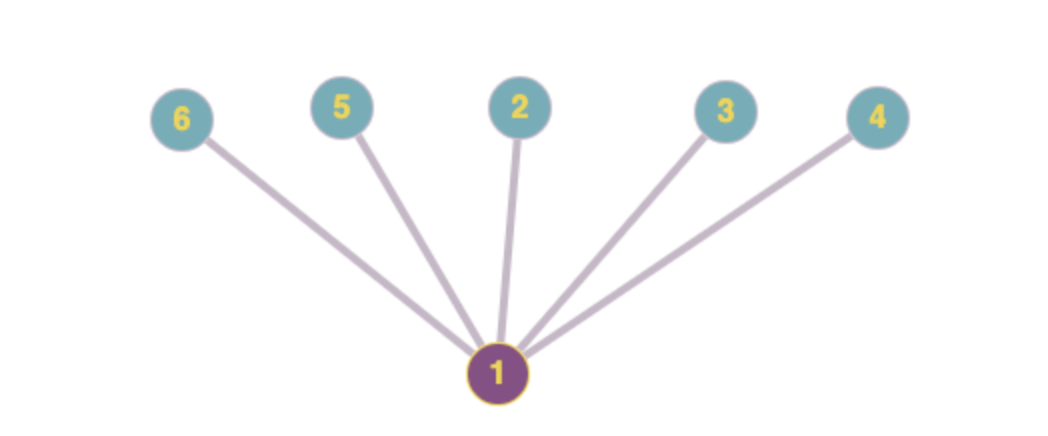}   
		\end{minipage}
	}
	\subfigure[bijective clique extension for $m=3$] 
	{
		\begin{minipage}{5.2cm}
			\centering      
			\includegraphics[scale=0.32]{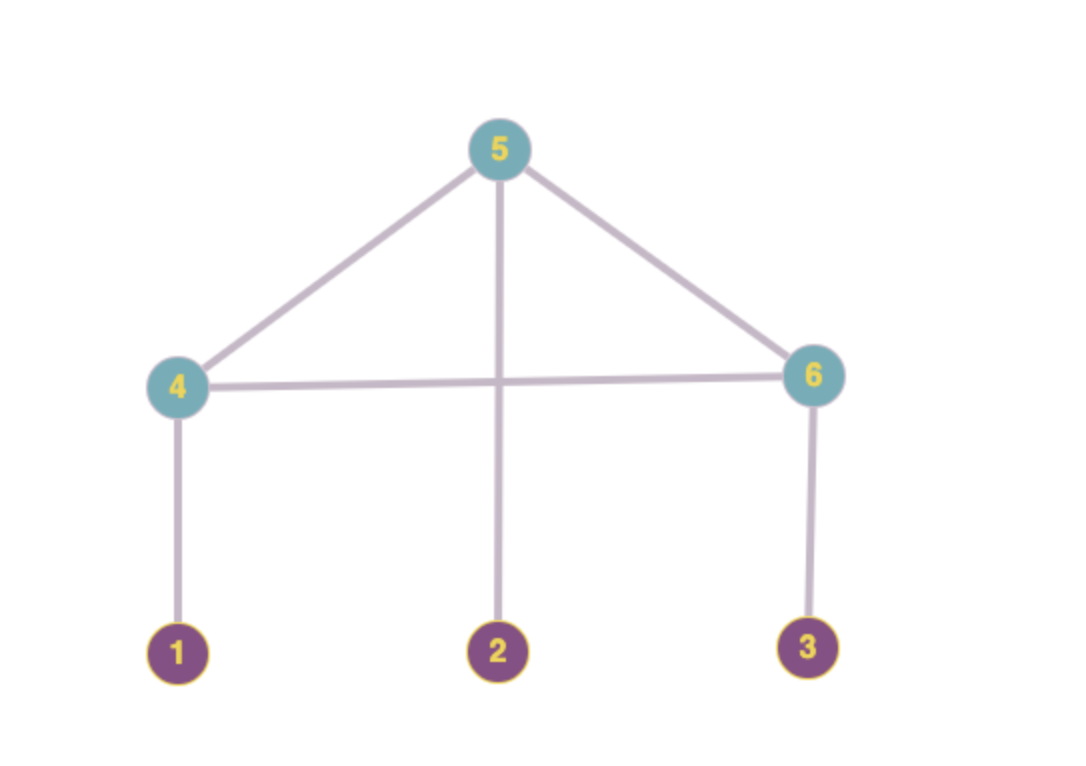}   
		\end{minipage}
	}
	\centering
	\subfigure[extension by path of length 5] 
	{
		\begin{minipage}{5.2cm}
			\centering      
			\includegraphics[scale=0.35]{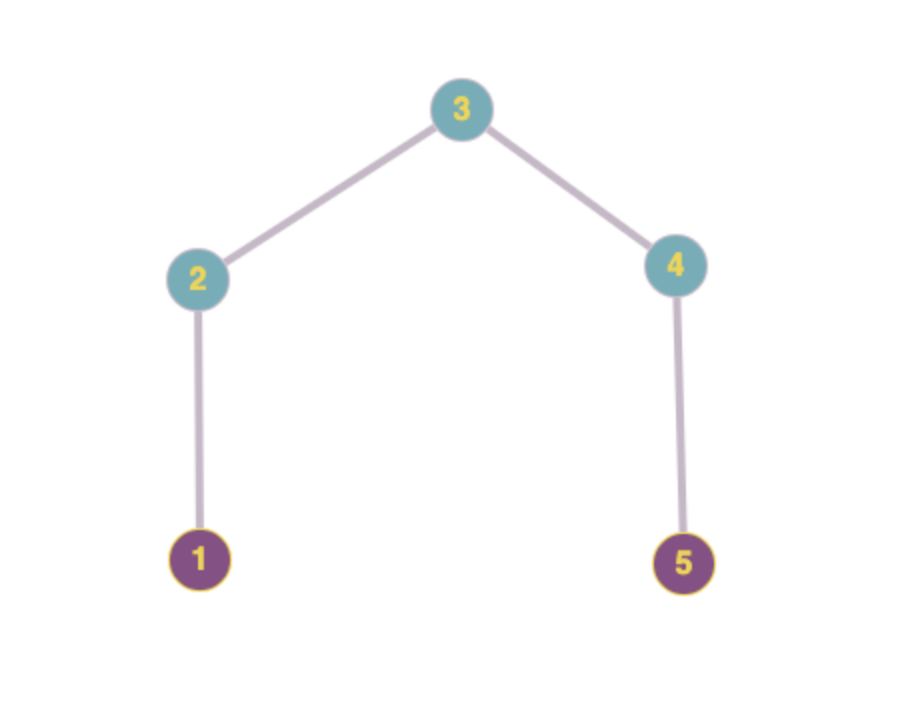}   
		\end{minipage}
	}
	\caption{symmetric rooted graphs, roots are in red} 
	\label{ext_examples}  
\end{figure}

In the current paper we consider $G(n, p=\mathrm{const})$ (in order to avoid hard technical details; however, at least some of our results can be generalised to a wider range of $p=p(n)$) and address the following general question. 
\begin{center}
Given a rooted graph $(R,H)$, what is the asymptotical distribution of the maximum of $X(T)$\\ over all possible choices of $r$-tuples $T$? 
\end{center}
More precisely, are there $a_n$ and $b_n$ such that $\frac{\max_T X(T)-a_n}{b_n}$ converges weakly to a non-generate distribution, and what is the limit distribution if this is the case? For convenience, we consider only {\it fully grounded} rooted graphs $(R,H)$, i.e. every root has at least 1 non-root neighbour. This assumption does not cause any loss in generality since clearly roots that are not adjacent to non-root vertices do not affect the maximum statistics we are looking at. In this paper we answer positively to the above question under certain conditions on $(R, H)$. More precisely, let us call $(R,H)$ \textit{symmetric}, if the set of root vertices $R$ can be divided into disjoint classes so that each non-root vertex is either not connected to the root set in $H$ or is connected to all vertices of exactly one class. So the expansion set $S(H)$ forms an arbitrary graph, and the only constraint is that the bipartite graph between $S(H)$ and $R$ is a disjoint union of complete bipartite subgraphs. Further in this section, we state the main result of our paper claiming a limit law for every symmetric extension. It generalises all the above mentioned results. Let us give various examples (see Fig.~\ref{ext_examples}) of symmetric rooted graphs including the three instances for which the limit law was known:

\renewcommand{\figurename}{Fig.}

a) $H$ is a single edge with a single root. Then $X_{(R, H)}(v) = \mathrm{deg}(v)$, the asymptotic distribution of the maximum degree was described in Theorem~\ref{th:max_deg}.

b) $H$ is a star graph with $k$ rays, all leaves are roots. In this case $X_{(R, H)}(v_1, \ldots, v_k) = \mathrm{deg}(v_1, \ldots , v_k)$, that denotes the number of common neighbours of vertices $v_1, \ldots , v_k$ in $G(n, p)$, the respective maximum was studied in~\cite{common}. In what follows, we denote by $\mathrm{deg}_G(U) $ and $N_G(U)$ the number of common neighbours and the set of common neighbours of vertices from the set $U$ in $G$ respectively.  We omit the subscript $G$, when the host graph $G$ is clear from the context.

c) $H$ is an $s$-clique with a single vertex being root. So $X_{(R, H)}(v)$ is the number of $s$-clicks that share $v$. The respective maximum was studied in \cite{main}.

Note that in the above three cases the bipartite graph between the set of roots and the expansion set is complete (i.e. there is a single class of roots), which appears to be crucial for the limit distribution to be from the Gumbel family. Let us give other two illustrative examples of symmetric extensions with several classes of roots:

d) $H$ consists of a set of roots and an expansion set of equal size $m$, the bipartite graph between them is a matching, and the expansion set induces an $m$-clique. We call such an extension \textit{a bijective ($m$-)clique extension}. Note that in this case there are $m$ classes of roots, each one consists of a single vertex.

e) $H$ is a simple path between two vertices $x_1, x_2$, the set of roots is $R = \{x_1, x_2\}$. There are exactly two classes of roots $\{x_1\}$ and $\{x_2\}$. Note that the respective maximum statistics is the maximum number of paths of a given length between a pair vertices.

As we will see later, the limit distributions of the maximum statistics related to the last two extensions do not belong to the Gumbel family.
\vspace{\baselineskip}

Let us now introduce the necessary notations and state the main result of our paper. Consider a symmetric fully grounded rooted graph $(R, H)$ with $h$ vertices and $f$ edges induced by the expansion set $S(H)$. Let its set of roots $R$ be divided into classes (in accordance with the definition of classes of roots of symmetric extensions) such that, for every $i\in[r]$, there are exactly $m_i$ classes of size $k_i$ (here, $k_1< \ldots < k_r$ are cardinalities of all the root classes that are presented in $H$). It turns out that the limiting distribution (but not the scaling constants) depends solely on the bipartite rooted subgraph of $H$ consisting of the same set of roots, vertices that are adjacent to at least one root in $H$ and edges between the roots and non-roots. This subgraph is defined by the vector $W(H):=((m_1, k_1)$, $(m_2, k_2), \ldots ,(m_r, k_r))$ as well as the vector of cardinalities of sets of vertices from the expansion set that are adjacent to all roots from a class (over all classes). Thus, to determine this subgraph completely, we consider $g_{ij}$, $i\in[r]$, $j\in[m_i]$, being the number of common neighbours of the $j$th root class of size $k_i$ in the expansion set. Without loss of generality we assume that $g_{i,1} \geq \ldots \geq g_{i, m_i}$ for every $i\in[r]$. Let us denote by $g_i:=\sum_{j=1}^{m_i} g_{i,j}$ the number of vertices adjacent to all roots from a certain class of size $k_i$, and by $g:=\sum_{i=1}^r g_i$ the total number of vertices adjacent to at least one root. Finally, let $s\geq 0$ be the number of vertices from the expansion set that are not adjacent to roots. Clearly, 

\begin{equation}\label{pattern_params_eq}
    |R| + g + s = h \hspace{5pt}.
\end{equation}

\begin{theorem} \it
    Within the above notations, define
    
    \begin{equation}\label{main_res_const}
        \begin{array}{l}
a_n=\frac{n^{s + g-1} p^f}{g_{1, 1}! g_{1, 2}! \ldots g_{r, m_r}!} \Big[n p^{\sum\limits_{i=1}^r k_i g_i} + \sqrt{2n \ln n} \times \\ \qquad\qquad\qquad\qquad\qquad \times \Big(\sum\limits_{i=1}^{r} g_i p^{k_i (g_{i} - 1)} \sqrt{k_i p^{k_i}(1-p^{k_i})}\Big(1 - \frac{\ln(k_i!)}{2k_i \ln n} - \frac{\ln[4\pi k_i \ln n]}{4k_i \ln n}\Big)\Big) \Big] ,\\
b_n=\frac{n^{s + g - 1} p^f}{g_{1, 1}! g_{1, 2}! \ldots g_{r, m_r}!} \sqrt{\frac{n}{2 \ln n}} p^{\sum\limits_{i=1}^r k_i g_i} .
        \end{array}
    \end{equation}
    Then
    
    \begin{equation}\label{main_res_lim}
        \frac{\max_T X(T) - a_n}{b_n} \xrightarrow[]{d} \sum\limits_{i=1}^r \sqrt{\frac{1-p^{k_i}}{k_i p^{k_i}}} \sum_{j=1}^{m_i} g_{i, j} \eta_{i, j} ,
    \end{equation}
where the vectors $\eta_i=(\eta_{ij},\,j\in[m_i])$ are mutually independent and have densities
$$p_{\eta_i}(x_{1}, x_{2}, \ldots ,x_{m_i}) = e^{-x_{1}} \cdot e^{-x_{2}} \cdot \ldots \cdot e^{-x_{m_i}} \cdot e^{-e^{-x_{m_i}}} \cdot I(x_{1} \geq x_{2} \geq \ldots \geq x_{m_i}) .$$
\label{th:main}
\end{theorem}
Let us now briefly discuss the methods of the proof. It seems natural that the maximum number of extensions is achieved at the set of roots whose classes have maximum number of common neighbours. For example, it turns out that the maximum number of paths of a given length is drawn between two vertices with the first and the second maximum degrees. In the same way, a pair of vertices with maximum number of common neighbours has maximum possible number of $k$-cliques inside its neighbourhood. This can be proven using a {\it conditional maximisation method} that we distill from~\cite{main} and develop and generalise in the present paper. In~\cite{main} in this way the limit distribution of the maximum number of $k$-cliques sharing a single vertex was studied. Let us briefly recall the main line of the proof. For every vertex $i$ of the random graph, consider its degree $\mathrm{deg}(i)$, and let $Y_i$ be the expected number of $k$-cliques containing $i$ conditioned on $\mathrm{deg}(i)$. The key argument that allows to transfer the limit distribution of $\max Y_i$ to the desired maximum number of $k$-cliques sharing a single vertex is

\begin{lemma}[M. Isaev, I. Rodionov, R. Zhang, M. Zhukovskii~\cite{main}]
\it
Let $X(n) \in \BR^d, d = d(n)$, be a sequence of random vectors, $a_n$ and $b_n$ --- two sequences of constants, and $F$ be a continuous cdf.

\noindent Let for any $x$ such that $0 < F(x) < 1$:
\begin{enumerate}

\item $ \prod_{i=1}^d {\sf Pr}( Y_i \leq a_n + b_n x) \rightarrow F(x)$,

\item ${\sf Pr}(\max_{i \in [d]} Y_i \leq a_n + b_n x) \rightarrow F(x)$,

\item for any fixed $\epsilon > 0$,
\begin{equation}\label{main_lemm}
{\sf Pr}(|X_i - Y_i| > \epsilon b_n) = o(1) {\sf Pr}(Y_i > a_n + b_n x) \text{\hspace{30pt}    uniformly over all $i \in [d]$}.
\end{equation}
\end{enumerate}

Then ${\sf Pr}(\max_{i \in [d]} X_i \leq a_n + b_nx) \rightarrow F(x)$ as well.
\label{lm:main_tech}
\end{lemma}

\vspace{\baselineskip}

In the present paper we generalise this techniques to symmetric rooted graphs with arbitrary root classes. This is possible since the conditional expectation is a monotone function of cardinalities of common neighbourhoods of root classes. For this reason, we find the limiting distribution of the vector of maximums $\Delta^j_{n, k_i}, i\in[r], j\in[m_i]$, where $\Delta^j_{n, k}$ is the $j$th maximum number of common neighbours of a $k$-set in $G(n, p)$. This generalises the main result of~\cite{common}. Note that, in particular, we show that whp the maximums are achieved at disjoint sets of roots ($m_1$ sets of size $k_1$, $m_2$ sets of size $k_2$, etc). Thus, this is possible to find explicitly the average number of $(R, H)$-extensions of these maximising sets of roots.

\vspace{\baselineskip}

Let us now apply Theorem~\ref{th:main} to rooted graphs described in a)-e). All these rooted graphs have $r=1$. 

Note that all the rooted graphs defined in a), b), c) have $m_1=1$ implying that the limit distribution belongs to the Gumbel family. In particular, consider a rooted graph with $k_1$ roots and $g$ pairwise adjacent non-roots, that are also adjacent to every root. This rooted graph generalises all rooted graphs from a), b), c). For the maximum number $\max_T X(T)$ of such extensions in $G(n, p)$ we get (we let $k=k_1$)

\begin{corol}\label{first_cor}
    Let $r=1$, $m_1=1$ and $s=0$. Let
$$
    \begin{array}{l}
a_n = \frac{(np^k)^{g-1} p^{g \choose 2}}{g!} \left[np^{k} + \sqrt{2n \ln n} g\sqrt{kp^k(1-p^k)}\left(1 - \frac{\ln(k!)}{2k \ln n} -  \frac{\ln[4\pi k \ln n]}{4 k\ln n}\right)\right] ,\\
b_n = \frac{n^{g-1} p^{{g \choose 2} + kg}}{(g-1)!} \sqrt{\frac{n(1-p^k)}{2k p^k \ln{n} }} .
    \end{array}
$$
    Then $\frac{\max_T X(T) - a_n}{b_n}\stackrel{d}\to\eta$, where $\eta$ has cdf $e^{-e^{-x}}$.
\end{corol}

Note that this number $\max_T X(T)$ is exactly the maximum number of $g$-cliques with at least $k$ common neighbours of their vertices. It is worth mentioning that this claim was announced in~\cite{common}, however its complete proof was not presented.

\vspace{\baselineskip}

Let us apply Theorem~\ref{th:main} to the case d). Here $W(H) = (m, 1)$, $v = 0, g_{1, j} = 1, g_{1} = g = m, f = {m \choose 2}$. By Theorem~\ref{th:main}, we get that the cdf of the limiting random variable equals

\begin{equation}\label{case1_ro}
    F(x) =  \int_{-\infty}^{x/m} \int_{t_m}^{(x-t_m)/(m-1)} \ldots \int_{t_2}^{x-t_m - \ldots - t_2} e^{-e^{-t_m}} e^{-t_m} e^{-t_{m-1}} \ldots e^{-t_1} dt_1 \ldots dt_m .
\end{equation}

 After accurate calculations, we can verify that its density function equals
 $$\rho(x) = e^{-x}\left(\frac{e^{-e^{-x/m}}}{m!} + \text{\\} P(x)\int_{-\infty}^{-e^{-x/m}} \frac{e^t}{t}dt\right)$$ for some polynomial $P$ since $F(x)$ can be represented as
\begin{equation}\label{int1}
\begin{array}{l}
    F(x)=\int\limits_{-\infty}^{x/m} e^{-e^{-t_m}} \frac{e^{-mt_m}}{(m-1)!}dt_m - e^{-x} \sum\limits_{i=2}^{m} \frac{1}{(i-1)!} I_i(x) \text{, where  } \\
    I_i(x) =  \int\limits_{-\infty}^{x/m} e^{-e^{-t_m}} \int\limits_{t_{m-1}}^{(x-t_m)/(m-1)} \ldots \int\limits_{t_{i+1}}^{(x-\sum\limits_{j=i+1}^{m} t_j)/i} dt_m \ldots dt_i.
\end{array}
\end{equation}

Note that $\mathrm{Ei}(y)= \int\limits_{-\infty}^{y} \frac{e^t}{t} dt$ is an exponential integral which is not an elementary function. Thus:

\begin{corol}\label{second_cor}
    Let $(R,H)$ be a rooted graph presented on Fig.~\ref{ext_examples}.d) with a clique of size $m\geq 2$. Let
     
$$
a_n = (np)^{m-1}p^{{m \choose 2}} \left[np + \sqrt{2n \ln n} m\sqrt{p(1-p)}\left(1 - \frac{\ln[4\pi \ln n]}{4\ln n}\right)\right],\quad
b_n = (np)^{m-1}p^{m \choose 2} \sqrt{\frac{np(1-p)}{2 \ln n}}.
$$
    
Then $\frac{\max_T X(T) - a_n}{b_n}\stackrel{d}\to\eta$, where $\eta$ has cdf described in $(\ref{int1})$.
\end{corol}

\vspace{\baselineskip}

Finally, we apply Theorem~\ref{th:main} to the case e), which corresponds to the maximum number of paths with $\ell>3$ edges between two vertices ($\ell=2, 3$ are special cases of Corollaries~\ref{first_cor}~and~\ref{second_cor} respectively). Here $W(H) = (2, 1)$, $v = \ell-3, g_{1, j} = 1, g_{1} = g = 2, f = {\ell-2}$. Note that the limit distribution is a particular case of (\ref{case1_ro}) with $m=2$ since, as we noted above, the limit distribution depends only on $W(H)$ and $(g_{ij})$, so its density equals

$$\rho(x) = \frac{d}{dx} \int_{-\infty}^{x/2} \int_{t_2}^{x-t_2} e^{-e^{-t_2}} e^{-t_2} e^{-t_1} dt_1 dt_2 = -e^{-x} \int_{-\infty}^{-e^{-x/2}} \frac{e^t}{t}dt.$$
Thus, we got the following result:

\begin{corol} Let $(R,H)$ be a rooted graph presented on Fig.~\ref{ext_examples}.e) with a path of length $\ell\geq 4$. Let
    $$
a_n = (np)^{\ell-2}p \left[np + 2\sqrt{2n \ln n p(1-p)}\left(1 - \frac{\ln[4\pi \ln n]}{4\ln n}\right)\right], \quad
b_n = (np)^{\ell-2}p \sqrt{\frac{np(1-p)}{2 \ln n}}.
$$
Then $\frac{\max_T X(T) - a_n}{b_n}\stackrel{d}\to\eta$, where $\eta$ has density $-e^{-x} \mathrm{Ei}(-e^{-x/2})$.
\end{corol}

So, indeed, the limit distributions of the maximum statistics from d) and e) does not belong to the Gumbel family.

\vspace{\baselineskip}

The rest of the paper is organised as follows. In Section~\ref{sc:pre} we recall and state several auxiliary claims about the random graph related to the binomial distribution that we use later in the proof. Section~\ref{sc:joint} is devoted to the joint limit distribution of scaled maximum numbers of common neighbours. The main result is proved in Section~\ref{sc:proof}. Section~\ref{sc:further} is devoted to a discussion of further questions.
\section{Preliminaries}
\label{sc:pre}

When working with maximum numbers of extensions, we frequently use asymptotical expressions for tails of binomial distribution from~\cite[Section 2.1]{common}, that follow from the de Moivre--Laplace limit theorem. In particular, the de Moivre--Laplace limit theorem immediately implies

\begin{claim}
Fix $\ell\in\mathbb{N}$ and $x > 0$. Consider arbitrary $\ell$ vertices $a_1, a_2, \ldots ,a_{\ell}$ in the random graph. Then

\begin{equation}
    {\sf Pr}\left(|\deg(a_1, \ldots ,a_{\ell}) - np^{\ell}| > \sqrt{2x np^{\ell}(1-p^{\ell})\ln n}\right) = \frac{1 + o(1)}{n^{x}\sqrt{\pi x \ln n}}.
\end{equation}
\label{cr:neigb}
\end{claim}

Let us denote for convenience $\Gamma_{\ell} = np^{\ell} + \sqrt{2{\ell} np^{\ell}(1-p^{\ell}) \ln n}$. By the union bound, the number of common neighbours of every set of $\ell$ vertices is at most $\Gamma_{\ell}$. Further in the work, in many places we restrict the probability space of graphs to only those graphs in which this property is satisfied for all $\ell \leq k$, where $k$ is a predefined fixed integer. We call this subspace $\mathcal{Q}_{n}$ (omitting the dependence of $k$ in the notation since it is always clear from the context), this narrowing would not affect convergences of probabilities to $0$ or $1$.

We also use the main result from~\cite{common} about the limit distribution of the maximum number of common neighbours.

\begin{theorem}[I. Rodionov, M. Zhukovskii~\cite{common}]
\it
Let $\Delta_{n, k}^m$ ($k, m \in \BN$) be the $m$-th highest number of common neighbours of $k$ vertices in $G(n, p)$,  where the maximum is taken over all possible $k$-tuples of distinct vertices. Let the probability of drawing an edge $p = p(n) \in (0, 1)$ be such that
$$
p^k \gg \frac{\ln^3 n}{n},\quad  1 - p \gg \sqrt{\frac{\ln \ln n}{n}}\quad\text{ as $\,\,n \rightarrow \infty$.}
$$

Let 
\begin{equation}\label{k_deg}
\begin{gathered}
a_{n, k} = np^k + \sqrt{2kp^k(1-p^k)n \ln n} \left( 1 - \frac{\ln(k!)}{2k \ln n} - \frac{\ln[4\pi k \ln n]}{4k \ln n}\right), \quad b_{n, k} = \sqrt{\frac{p^k(1-p^k)n}{2k \ln n}}.
\end{gathered}
\end{equation}

Then $\frac{\Delta_{n,k}^m - a_{n, k}}{b_{n, k}}$ converges in distribution to a random variable with cdf $e^{-e^{-x}} \sum \limits_{j=0}^{m-1} \frac{e^{-jx}}{j!}$.
\label{th:max_com_neigb}
\end{theorem}

We also use the asymptotics of the probability that a fixed $k$-set $U$ has more than $a_{n,k} + xb_{n,k}$ common neighbours. Denoting this event by $B_U(x)$, using the de Moivre--Laplace limit theorem, it is easy to see (the full proof can be found in~\cite[Section 2.1]{common}) that

\begin{equation}\label{conv_exp}
\begin{gathered}
    {\sf Pr}(B_U(x)) \sim \frac{k!}{n^k} e^{-x} \text{ as } n \rightarrow{}{} \infty.
\end{gathered}
\end{equation}

In Appendix, we prove the useful technical lemma which is stated below. It claims that the maximum numbers of common neighbours are achieved at non-overlapping sets. We use this lemma to show that the maximum number of extensions is achieved at those disjoint root classes that, in turn, admit maximum numbers of respective subextensions by common neighbours. 

\begin{lemma}
\it Let $m_i, k_i \in \mathbb{N}$, $i\in[r], r\in\mathbb{N}$, and all $k_i$ be distinct. Let $U_{i,j}$, $i\in[r], j\in[m_i]$, be $k_i$-sets such that cardinalities of their common neighborhoods are maximum, i.e. for every $i\in[r]$ $\deg(U_{i,1}) \geq \ldots\geq \deg(U_{i, m_i})$ are cardinalities of $m_i$ biggest common neighborhoods among all $k_i$-sets. Then whp all $U_{i,j}$ are disjoint.
\label{lm:intersection}
\end{lemma}

We move the proof to Appendix B since it is actually a generalisation of a particular case of this result proven (implicitly) in~\cite{common}, and we use exactly the same proof strategy.
\section{ Joint distribution of maxima}
\label{sc:joint}

The limit distribution of the scaled maximum number of extensions in Theorem~\ref{th:main} is in fact entirely determined by the joint distribution of the maximum numbers of common neighbours of sets of vertices of respective sizes, which is studied in this section. In the first subsection, we find the joint distribution of $\Delta_{n,k_i}:=\Delta^1_{n,k_i}$, $i\in[r]$, --- maxima cardinalities of common neighborhoods of $k_i$ vertices for distinct $k_1,\ldots,k_r$. 
In the second subsection, using this result, we find the limit joint distribution of the first $m_i$ largest numbers of common neighbours of $k_i$ vertices, $i \in [r]$.

\subsection{Maximum neighborhoods}
It is shown here that the scaled maximum numbers of common neighbours are almost independent. More precisely, the following generalisation of Theorem~\ref{th:max_com_neigb} (for constant $p$) is proved:

\begin{claim}
\it
Let some $x_1, x_2, \ldots , x_r \in \BR$ be fixed. Then
\begin{center}
${\sf Pr}\left(\frac{\Delta_{n, k_1} - a_{n, k_1}}{b_{n, k_1}} \leq x_1, \frac{\Delta_{n, k_2} - a_{n, k_2}}{b_{n, k_2}} \leq x_2, \ldots, \frac{\Delta_{n, k_r} - a_{n, k_r}}{b_{n, k_r}} \leq x_r\right) \rightarrow e^{-e^{-x_1}} \cdot e^{-e^{-x_2}} \cdot \ldots \cdot e^{-e^{-x_r}}$ as $n\to \infty$,
\end{center}
where constants $a_{n, k_i}, b_{n, k_i}$ are defined in (\ref{k_deg}).
\label{st:maximum}
\end{claim}

Denote by $X_i = X_i(x_i)$, $i \in [r]$, the number of sets of $k_i$ vertices that have a ``large'' number of common neighbours, namely, more than $a_{n, k_i} + b_{n, k_i }x_i$. Then our goal is to bound ${\sf Pr}(X_1 = 0, X_2 = 0, \ldots, X_r = 0)$.\\

\textbf{Lower bound} 
\begin{equation}
\nonumber
    {\sf Pr}(X_1 = 0, X_2 = 0, \ldots, X_r = 0) \geq {\sf Pr}(X_1 = 0) {\sf Pr}(X_2 = 0) \ldots {\sf Pr}(X_r = 0)
\end{equation}

\noindent is a consequence of \cite[Theorem 6.3.3]{prob_method} --- an application of the well-known FKG-inequality~\cite[Theorem 6.2.1]{prob_method}. Indeed, the properties of the absence of sets with a large number of common neighbours are \textit{decreasing} functions of the edges of the random $G(n, p)$. The limit of the right-hand side of this bound coincides with the limit distribution in Claim~\ref{st:maximum} due to Theorem~\ref{th:max_com_neigb}.\\

\textbf{Upper bound} is in fact similar to the proof of \cite[Lemma 1]{common} and follows almost directly from~\cite[Lemma 3.1]{main}. Let us recall the requirements and the statement of this lemma. 

Let us denote by $T$ the set of all subsets of vertices in $G(n, p)$ of one of the sizes $k_1, k_2, \ldots, k_r$. We consider two families of events: $\{B_U\}$ and $\{\tilde B_U\} = \{B_U \cap \{G \in \mathcal{Q}_{n}\} \},$ where $U = \{u_1, \ldots ,u_{k_i}\}$, $i\in[r]$, is an arbitrary set in $T$. Note that $x_i$ is substituted into the definition of $B_U=B_U(x_i)$ according to the size of $U$. Thus our aim is to bound ${\sf Pr}( \bigcap_{U \in T} \overline{B_U}) \leq {\sf Pr}( \bigcap_{U \in T} \overline{\tilde{B_U}})$. To do this, we use the following key lemma from ~\cite[Lemma 3.1]{main}.

\begin{lemma}[M. Isaev, I. Rodionov, R. Zhang, M. Zhukovskii~\cite{main}]
	Let $(A_i)_{i \in [d]}$ be the set of events with non-zero probabilities. If sets $(D_i \subset [d] \backslash \{i\})_{i\in[d]}$ satisfy
	$${\sf Pr} \left(  \bigcup_{j \in [i-1] \backslash D_i} A_j | A_i \right) - {\sf Pr} \left(  \bigcup_{j \in [i-1] \backslash D_i} A_j \right)  \leq \varphi ,$$
	for some $\varphi \geq 0$ and all $i\in[d]$, then
	\begin{equation}
        {\sf Pr}\left(\bigcap_{i \in [d]} \overline{A_i}\right) \leq \prod_{i\in[d]} {\sf Pr}(\overline{A_i}) + \varphi \left( 1 - \prod_{i \in [d]} {\sf Pr} ( \overline{A_i}) \right) + \Delta ,
        \label{lm:good}
	\end{equation}
	where $\Delta = \Delta(A, D) = \sum\limits_{i\in[d]} {\sf Pr}\left( A_i \cap \bigcup\limits_{j \in [i-1]\cup D_i} A_j\right) \prod\limits_{\ell \in [d] \backslash [i]}{\sf Pr} (\overline{A_{\ell}})$.
	\label{lm:formulation:good}
\end{lemma}

It is useful to choose $D_i$ to be the set of all $j\neq i$ so that $A_j$ strongly depends on $A_i$. We order all $U\in T$, and let $A_i=\tilde B_U$ for the $i$th set $U$. We also let $j\in D_i$ whenever the $j$th set of $T$ has a non-empty intersection with the $i$th set from $T$. Then
$$\prod {\sf Pr} (\overline{A_i})=\prod_{U\in T} {\sf Pr}(\overline{\tilde{B}_U}) = \exp \left[ \sum\limits_{U\in T} \ln (1 - {\sf Pr}(\tilde{B}_U)) \right] = \exp \left[ \sum\limits_{i\in [r]} -\lambda_{k_i} + o(1) \right], $$
where $\lambda_{k} = \sum\limits_{U \subset [n], |U| = k} {\sf Pr} (\tilde B_U)$. In \cite[Section 2.3.1]{common} it is proved that $\lambda_k \sim e^{-x_k}$ as $n \rightarrow \infty$. Thus, it suffices to verify that $\Delta=o(1)$ and $\varphi=o(1)$.

Let us first prove that $\Delta=o(1)$. In the proof of Lemma~\ref{lm:intersection} it is shown that for arbitrary $i, j \in [r]$ and an arbitraty $C \in \mathbb{R}$ 
$$\sum\limits_{U \cap V \neq \varnothing, |U| = k_i, |V| = k_j, U \neq V} {\sf Pr}\left(\deg(U) > a_{n, k_i} + C \sqrt{\frac{n}{\ln n}},\,\, \deg(V) > a_{n, k_j} + C \sqrt{\frac{n}{\ln n}},\,\, G \in \mathcal{Q}_{n}\right) \rightarrow 0 .$$
Choose $C$ sufficiently small and get
$$
 \Delta\leq\sum_{U\in T,V\in T:\,V\cap U\neq\varnothing}\Pr(\tilde B_U\cap\tilde B_V)=o(1) .
$$
In remains to prove that $\phi=o(1)$. For every $U\in T$
\begin{multline*}
{\sf Pr} \left(  \bigcup_{V \cap U = \varnothing} \left. \tilde{B}_V \right| \tilde{B}_U \right) - {\sf Pr} \left(  \bigcup_{V \cap U = \varnothing} \tilde{B}_V \right) \leq \\
\leq{\sf Pr} \left(  \bigcup_{V \cap U = \varnothing} \deg_{G \backslash U}(V) > a_{n,k_i} + x_i b_{n, k_i} - |U| \right) - {\sf Pr} \left(  \bigcup_{V \cap U = \varnothing} \tilde{B}_V \right) ,
\end{multline*}
where $k_i=k_i(V)=|V|$ and $x_i=x_i(V)$ is defined accordingly. So due to the union bound and the de Moivre--Laplace limit theorem we get
$${\sf Pr} \left(  \bigcup_{V \cap U = \varnothing} \left. \tilde{B}_V  \right| \tilde{B}_U \right) - {\sf Pr} \left(  \bigcup_{V \cap U = \varnothing} \tilde{B}_V \right) \leq \sum\limits_{V \subset T} {\sf Pr} (\deg(V) \in [-k_r, 0] + a_{n,k_i} + x_i b_{n, k_i} ) \rightarrow 0$$
uniformly over $i \in [d]$, implying that $\varphi = o(1)$ and completing the proof.
\subsection{First $m_i$ maxima}

For $i\in[r]$ and $j \in[m_i]$, let $\xi_{i, j}$ be the centered and normalised $j$-th maximum number of common neighbours of $k_i$ vertices in $G(n, p)$ with the scaling constants defined in (\ref{k_deg}), i.e. $$\xi_{i, j} = \frac{\Delta_{n, k_i}^ {j} - a_{n, k_i}}{b_{n, k_i}} .$$ The purpose of this section is to find the limiting distribution of the random vector $\xi$ comprising all $s = \sum\limits_{i=1}^r m_i$ random variables $\xi_{i, j}$, $i\in[r]$, $j \in [m_i]$.

For $x\in\mathbb{R}^s$ we will denote its coordinates by $x_{i, j}$, $i\in[r]$, $j \in [m_i]$, for convenience. Clearly, it is sufficient to study the distribution of $\xi$ on the set $Y = \{x \in \BR^s : \forall i \in [r]$ $x_{i, m_i} \leq x_{i , m_i - 1} \leq \ldots \leq x_{i, 1}\}$, since from the definition $\xi_{i, m_i} \leq \xi_{i, m_i-1} \ldots \leq \xi_{i, 1}$ for every $i \in [r]$. Fix $x\in\mathbb{R}^s$. For $i\in[r]$, set $A(i) = \{\xi_{i, 1} \leq x_{i, 1}, \xi_{i, 2} \leq x_{i, 2}, \ldots , \xi_{i, m_i} \leq x_{i, m_i}\}$.

For $i\in[r]$, $t\in[m_i]$ and $1\leq\ell_1\leq\ell_2\leq\ldots\leq\ell_{t-1}\leq m_i$, define 
$$A(i; \ell_1, \ldots, \ell_{t-1})=\{\xi_{i, 1} \in [x_{i, \ell_1}, x_{i, \ell_1+1}], \ldots, \xi_{i, t-1} \in [x_{i, \ell_{t-1}}, x_{i, \ell_{t-1}+1}], \xi_{i, t} \leq x_{i, m_i}\}$$ 
--- the event, saying that each $\xi_{ij}$ (but the smallest one) is between two consecutive coordinates of $x$. Clearly, $A(i)$ is the disjoint union of all possible $A(i; \ell_1, \ldots, \ell_{t-1})$. So, in order to find the distribution of $\xi$ it is sufficient to find it on all Cartesian products of events $A(i; \ell_1, \ldots, \ell_{t-1})$ over $i\in[r]$. As we will see later, in order to compute density of the limit distribution of $\xi$, it is sufficient to find the measure of a one ``simple brick'' $D = D_1 \times \ldots \times D_{r}$, where:
\begin{equation}
\label{interesting_event}
\begin{gathered}
    D_i = \{\xi_{i, 1} \in [x_{i, 2}, x_{i, 1}],\xi_{i, 2} \in [x_{i, 3}, x_{i, 2}], \ldots ,
    \xi_{i, m_i-1} \in [x_{i, m_i}, x_{i, m_i-1}], \xi_{i, m_i} \leq x_{i, m_i}\} .
\end{gathered}
\end{equation}

Let us also restrict the probability space only to those graphs in which the first $m_i$ maxima numbers of common neighbours of $k_i$-sets are reached at non-overlapping sets over all $i\in[r]$. We denote this event as $DisjRoots$. From Lemma~\ref{lm:intersection} whp $DisjRoots$ happens, so the limit of ${\sf Pr}(D_1 \times D_2 \times \ldots \times D_r)$ is the same as the probability limit of $D' = D_1 \times D_2 \times \ldots \times D_r \cap DisjRoots$.

Now we consider the set of disjoint events $D'(w), w\in W$, where $W = (U_{ij}, i\in[r], j\in[m_i-1])$ --- the set of all tuples of disjoint sets $U_{i,j}$ of size $k_i$, and
\begin{equation}
\nonumber
\label{each_event}
\begin{gathered}
D'(w) =  \bigcap\limits_{i=1}^r \left\{\forall j\in[m_i-1] \text{ } \frac{\deg(U_{i, j}) - a_{n, k_i}}{b_{n, k_i}} \in [x_{i, j+1}, x_{i, j}], \text{ } \frac{\max\limits_{V_i \in G/U, |V_i| = k_i}{} \deg(V_i) - a_{n, k_i}}{b_{n, k_i}} \leq x_{i, m_i} \right\},
\end{gathered}
\end{equation}
where $U = \bigsqcup\limits_{ i\in[r], j\in[m_i-1]} U_{i,j}$. It is obvious that
$$\sum\limits_w {\sf Pr}(D'(w)) - {\sf}{\sf Pr}(DisjRoots) \leq {\sf Pr}(D') \leq \sum\limits_w {\sf Pr}(D'(w)),$$
so it is enough to estimate the sum of ${\sf Pr}(D'(w))$ over $w\in W$. The total number of vectors in $W$ is 
\begin{equation}\label{W_size}
\begin{gathered}
    |W| = {n \choose k_1} \cdot {n - k_1 \choose k_1} \cdot \ldots \cdot {n - (m_1 - 2)k_1 \choose k_1} \cdot {n - (m_1 - 1)k_1 \choose k_2} \cdot \ldots \cdot \\ \cdot {n - (m_1 - 1)k_1 - (m_2-2)k_2 \choose k_2} \cdot  \ldots \cdot {n - \sum (m_i-1)k_i + k_r \choose k_r} = \frac{n^{|U| }(1+o(1))}{(k_1!)^{m_1-1} (k_2!)^{m_2-1} \ldots (k_r!)^{m_r-1}} .
\end{gathered}
\end{equation}
Let us order pairs $(i, j)$ lexicographically. Denote $G_{ij} = G/\bigcup \limits_{(i',j') < (i,j)}U_{i',j'}$. Then we have for each $w \in W$:
\begin{equation}
\nonumber
\begin{gathered}
    D'(w) = \left\{\forall i\in[r] \text{  }  \forall j\in[m_i-1] \hspace{12pt} \frac{\deg_{G_{ij}}(U_{i, j}) - a_{n, k_i} + \epsilon_{i, j}}{b_{n, k_i}} \in [x_{i, j+1}, x_{i, j}], \right. \\
    \left. \forall i\in[r] \hspace{12pt} \frac{\max\limits_{V_i \in {G/U \choose k_i}}{} \deg_{G/U}(V_i) - a_{n, k_i} + \epsilon_{i}}{b_{n, k_i}} \leq x_{i, m_i} \right\},
\end{gathered}
\end{equation}
where $\epsilon_{i,j}$ and $\epsilon_i$ are random variable equal to the number of common neighbours of $U_{i,j}$ and  $V_i$ respectively among the union of the previous ones in the our enumeration $\{U_{i, j}\}$. It is clear that for all $i\in[r], j\in[m_i]$, $\epsilon_i, \epsilon_{i, j} < |R|=\mathrm{const}$. Using this, and the  consequence of the De Moivre-Laplace theorem (\ref{conv_exp}) we get that the probability limit is
\begin{equation} \nonumber
    \begin{gathered}
        \lim\limits_{n\to\infty} {\sf Pr}(D'(w)) = \prod\limits_{i=1}^r \left( \frac{(k_i)! (e^{-x_{i,2}} - e^{-x_{i, 1}})}{n^{k_i}} \times \ldots \times \frac{(k_i)!(e^{-x_{i, m_i}} - e^{-x_{i, m_i-1}})}{n^{k_i}} \right) \times \\
        \times \lim\limits_{n\to\infty} {\sf Pr}\left(\frac{\max\limits_{V_1 \in {G/U \choose k_1}}{} \deg_{G/U}(V_1) - a_{n, k_1}}{b_{n, k_1}} \leq x_{1, m_1}, \ldots, \frac{\max\limits_{V_r \in {G/U \choose k_r}}{} \deg_{G/U}(V_r) - a_{n, k_r} }{b_{n, k_r}} \leq x_{r, m_r}\right).
    \end{gathered}
\end{equation}

Using the probability limit for the last factor from Claim~\ref{st:maximum} and the asymptotics on $|W|$ (\ref{W_size}), we get
\begin{align*} 
    \lim\limits_{n\to\infty} {\sf Pr}(D) &= \lim\limits_{n\to\infty} {\sf Pr}(D') = \sum\limits_{w \in W}  \lim\limits_{n\to\infty} {\sf Pr}(D'(w)) = \\
    &= \prod\limits_{i=1}^r \left( (e^{-x_{i, 2}} - e^{-x_{i, 1}}) \cdot (e^{-x_{i, 3}} - e^{-x_{i, 2}}) \cdot \ldots  \cdot (e^{-x_{i, m_i}} - e^{-x_{i, m_i - 1}}) \right) \cdot \prod\limits_{i=1}^r e^{-e^{-x_{i, m_i}}} := F(x) \hspace{4pt}.
\end{align*}
We denote by $\mathcal{A}$ the set of all Cartesian products of $A(i; \ell_1, \ldots, \ell_{t-1})$ over $i \in [r]$. In the same way as above, it is easy to see that the limit probability of the $j$-th set $A_j = A(1; \ell^1_{1}, \ldots, \ell_{t_1-1}^1)\times \ldots \times A(r; \ell_{1}^r, \ldots, \ell_{t_r-1}^r) \in \mathcal{A}$ is
\begin{equation}
\nonumber
\begin{gathered}
\label{each_term}
T_j(x) := \prod\limits_{i=1}^r (e^{-x_{i, \ell_{1}^i+1}} - e^{-x_{i, \ell_{1}^i}}) \cdot (e^{-x_{i, \ell_{2}^i+1}} - e^{-x_{i, \ell_{2}^1}}) \cdot \ldots \cdot (e^{-x_{i, \ell_{(t_1-1)}^i+1}} - e^{-x_{i, \ell_{(t_1-1)}^i}}) \cdot \prod\limits_{i=1}^r e^{-e^{-x_{i, m_i}}} .
\end{gathered}
\end{equation}
 
 It is easy to see that the density of limit distribution of $\xi$ equals 
 \begin{align*}
 p(x_1, \ldots, x_s) &= \frac{\partial^s}{\partial x_1 \ldots \partial x_s} \sum \limits_{j=1}^{|\mathcal{A}|} T_j(x_1, \ldots, x_s) = \frac{\partial^s}{\partial x_1 \ldots \partial x_s} F(x_1, \ldots, x_s) =  \\ 
 &= \frac{\partial^s}{\partial x_1 \ldots \partial x_s}  \prod\limits_{i=1}^r (e^{-x_{i, 2}} - e^{-x_{1, 1}}) (e^{-x_{i, 3}} - e^{-x_{1, 2}}) \cdot \ldots \cdot (e^{-x_{i, m_i}} - e^{-x_{i, m_i - 1}}) \cdot e^{-e^{-x_{i, m_i}}} .
 \end{align*}
Expanding all brackets and differentiating, we obtain

\begin{claim}
\it
$\xi$ converges in distribution to a random vector with an absolutely continuous distribution with pdf $p(x_1, \ldots, x_s) = \prod\limits_{i=1}^r p_{i}(x_{i, 1}, x_{i, 2}, \ldots x_{i, m_i})$, where each
$$p_{i}(x_{1}, x_{2}, \ldots ,x_{m_i}) = e^{-x_{1}} \cdot e^{-x_{2}} \cdot \ldots \cdot e^{-x_{m_i}} \cdot e^{-e^{-x_{m_i}}} \cdot I(x_{1} \geq x_{2} \geq \ldots \geq x_{m_i}) .$$
\label{st:neighb_dist}
\end{claim}
Note that Theorem~\ref{th:max_deg} and Theorem~\ref{th:max_com_neigb} are particular cases of Claim~\ref{st:neighb_dist} for constant $p$.
\section{Proof of the main result}
\label{sc:proof}

In this section we prove the main result of the paper, Theorem~\ref{th:main}, by implementing the conditional maximisation method described in Introduction.  Let us consider in $G(n,p)$ an arbitrary ordered set of vertices $T$ of cardinality $|R|$ and its partition into root classes $A_{ij}, i\in[r], j\in[m_i]$. Let $Y(T)$ be the number of $(R, H)$-extensions conditioned on numbers of common neighbours for all root classes $A_{ij}$. Thus
\begin{equation}
\label{cond_ext} 
Y(T) = {\sf E} ( \left. X_{(R, H)}(T) \right|  \deg_G(A_{ij}), i\in[r], j\in[m_i]) .
\end{equation}
The general idea is to find the limit distribution of a scaled $\max\limits_T Y(T)$ and then prove that the maximum number of extensions $\max\limits_T X(T)$ is not much different from it and so converges to the same distribution. It is worth noting that we can not do the same as in~\cite{main} and directly apply Lemma~\ref{lm:main_tech} since the first condition is not satisfied in our settings: the product of probabilities does not converge to the limit distribution of maxima. However, we state a more general lemma, which is sufficient for our purposes:

\begin{lemma}
\it
Let $X=X(n) \in \BR^d, d = d(n)$, be a sequence of random vectors. Let $a_n$ and $b_n$ be two sequences of real constants, and let $F$ be a continuous cdf. Let, for any $x\in\mathbb{R}$ such that $0 < F(x) < 1$,
\begin{enumerate}

\item ${\sf Pr}(\max_{i \in [d]} Y_i \leq a_n + b_n x) \rightarrow F(x)$,

\item for any fixed $\varepsilon > 0$,
\begin{equation}\label{our_lemm}
\sum_{i=1}^d{\sf Pr}(|X_i - Y_i| > \varepsilon b_n) = o(1) \text{\hspace{2pt}  .}
\end{equation}
\end{enumerate}
Then ${\sf Pr}(\max_{i \in [d]} X_i \leq a_n + b_nx) \rightarrow F(x)$ for all $x\in\mathbb{R}$.
\label{lm:our_tech}
\end{lemma}

The proof of this lemma is similar to the proof of Lemma~\ref{lm:main_tech}; it can be found in Appendix A. We verify the first requirement in Lemma~\ref{lm:our_tech} with cdf defined in (\ref{main_res_lim}) in Section 4.1. The second condition is verified using Janson inequality and a similar (but weaker) upper tail bound in Section 4.2 completing the proof of Theorem~\ref{th:main}.

\subsection{Convergence of the expected conditional number of extensions}

Here we will havily rely on Claim~\ref{st:neighb_dist}.

Consider an arbitrary set of vertices $T$ of size $|R|$  and its partition in accordance with $W(H)$:

$$T = \bigsqcup\limits_{i=1}^r A_i, \hspace{30pt} A_i = \bigsqcup\limits_{j=1}^{m_i} A_{i,j}, \hspace{10pt} \text{where   } |A_{i,j}| = k_i .$$

Then for $Y(T)$ defined in (\ref{cond_ext}) we have:
$$Y(T) = p^f (n-h+s) \cdot (n-h+s-1) \cdot \ldots \cdot (n-h+1) \cdot  {\sf E}(\left. S(T) \right| \deg(A_{i,j}), i\in[r], j\in[m_i]) ,$$
where $S(T)$ is the number of $(R, H')$-extensions of $T$ in $G(n, p)$, and $H'$ is obtained from $H$ by deleting all non-root vertices that are not adjacent to roots and also all edges between all the remaining non-root vertices.

Let us estimate the conditional expectation of $S(T)$. From the definition of \textit{symmetric extensions}, each vertex of this ``first'' level in $H$ is connected to exactly one of the sets of roots corresponding to $A_{i,j}$ in $G(n,p)$. Note that, if $U$ is the set of all common neighbours of $A_{1, 1}$ in $G(n, p)$, then it may happen that some other $A_{i,j}$ has common neighbours in $U$ or that some roots from $T$ belong to $U$. Then obviously
$$\prod\limits_{i\in[r],j\in[m_i]} {\deg(A_{i,j}) - |R| - g \choose g_{i,j}} \leq {\sf E} \left(\left.S(T) \right| \deg(A_{i,j}), i\in[r], j\in[m_i])\right) \leq \prod\limits_{i\in[r],j\in[m_i]} {\deg(A_{i,j}) \choose g_{i,j}} .$$
Thus, assuming that all $\deg(A_{i,j})\to\infty$ as $n\to\infty$, we get that
\begin{equation}
Y(T)=p^f n^s \prod\limits_{i\in[r],j\in[m_i]} \frac{\deg(A_{i,j})^{g_{i,j}}}{g_{i,j}!}\left(1+O\left(\frac{1}{\deg(A_{i,j})}\right)\right) .
\label{both_bounds}
\end{equation}

Denote $\psi_{i, j} = \frac{\deg(A_{i, j}) - a_{n, k_i}}{b_{n, k_i}}$ with constants $a_{n, k}, b_{n, k}$ defined in (\ref{k_deg}). Note that for every $i\in[r]$, the first $m_i$ maxima of $\psi_{i,j}$ over $A_{i,j}$ equal $\xi_{i, 1} \geq \xi_{i,2} \ldots \geq \xi_{i,m_i}$, where $\xi_{i,j}$ are defined in Section 3.2. Since $a_{n, k} \sim n, b_{n, k} \sim \sqrt{\frac{n}{\ln n}}$, and whp $\psi_{i,j} = O(\ln n)$ (we further restrict the space of graphs to those in which this condition is satisfied, the convergence of probabilities does not change), then whp
$$Y(T) = a(n) + b(n) + o\left(n^{s+g-1}\sqrt{\frac{n}{\ln n}}\right),$$
where
\begin{equation}
\nonumber
\begin{split}
    a(n) & = \frac{p^f n^s}{\prod\limits_{i\in[r],j\in[m_i]} g_{i,j}!} \prod\limits_{i=1}^r a_{n, k_i}^{\sum\limits_{j=1}^{m_i} g_{ij}} \sim \frac{p^f n^{g+s-1}}{\prod\limits_{i\in[r],j\in[m_i]} g_{i,j}!} \left( np^{\sum\limits_{i=1}^r k_i \sum\limits_{j=1}^{m_i} g_{i,j}} + \right.\\ 
    & +  \left. \sqrt{2n \ln n} \left(\sum\limits_{i=1}^{r} \left(\sum\limits_{j=1}^{m_i} g_{ij}\right) p^{k_i \left(\sum\limits_{j=1}^{m_i} g_{i,j} - 1\right)} \sqrt{k_ip^{k_i}(1-p^{k_i})}\left(1 - \frac{\ln(k_i!)}{2k_i \ln n} - \frac{\ln[4\pi k_i \ln n]}{4k_i \ln n}\right)\right)\right) = a_n , \\
b(n) &= \frac{p^fn^s}{\prod\limits_{i\in[r],j\in[m_i]} g_{i,j}!} \sum\limits_{i=1}^r \left[ a_{n, k_i}^{\sum\limits_{j=1}^{m_i} g_{i,j} - 1} \left( \prod\limits_{i' \neq i}^r a_{n, k_{i'}}^{\sum\limits_{j=1}^{m_{i'}} g_{i',j}} \right) b_{n, k_i} \left(\sum\limits_{j=1}^{m_i} g_{i,j} \psi_{i, j}\right)\right] \\
& \sim \frac{p^fn^{s+g-1}}{\prod\limits_{i\in[r],j\in[m_i]} g_{i,j}!} \sqrt{\frac{n}{2 \ln n}} p^{\sum\limits_{i=1}^r k_i \sum\limits_{j=1}^{m_i} g_{i,j}} \left( \sum\limits_{i=1}^r \sqrt{\frac{1-p^{k_i}}{k_i p^{k_i}}} \sum_{j=1}^{m_i} g_{i, j} \psi_{i, j}\right) \\
& = b_n \left(\sum\limits_{i=1}^r \sqrt{\frac{1-p^{k_i}}{k_i p^{k_i}}} \sum_{j=1}^{m_i} g_{i, j} \psi_{i, j}\right) .
\end{split}
\end{equation}

By Claim~\ref{st:neighb_dist} and Slutsky's theorem,
\begin{equation}
\max_T \frac{Y(T) - a_n}{b_n} =  \sum\limits_{i=1}^r \sqrt{\frac{1-p^{k_i}}{k_i p^{k_i}}} \sum_{j=1}^{m_i} g_{i, j} \xi_{i, j} + o_{{\sf P}}(1) \xrightarrow[]{d} \eta ,
\label{result_convergence}
\end{equation}
where $\eta$ has cdf defined in (\ref{main_res_lim}). Note that the equality in (\ref{result_convergence}) holds true due to the descending order of $g_{i,j}$ for each fixed $i$ since $\xi_{i,1}\geq\ldots\geq\xi_{i, m_i}$. It is also worth noting that whp the maximum of $Y(T)$ coincides with the point-wise maximum (i.e. is achieved at $A_{i,j}$ that have maximum numbers of common neighbours). Finally, Lemma~\ref{lm:intersection} together with (\ref{result_convergence}) imply the first requirement in Lemma~\ref{lm:our_tech}.

\begin{remark}
The pdf of $\eta$ could be found explicitly due to Claim~\ref{st:neighb_dist}. Note that in the case $r=1$, we may divide both parts of $(\ref{result_convergence})$ by $\sqrt{\frac{1-p^{k_1}}{k_1 p^{k_1}}}$ avoiding the dependency of the limit distribution of $p$.

\end{remark}

\subsection{Deviation from the expected conditional number of extensions}
Here, using Janson-type correlation inequalities, we check the condition (\ref{our_lemm}):
$$\sum\limits_{T} {\sf Pr}(|X(T) - Y(T)| > \varepsilon b_n) = o(1) \text{\hspace{1pt}  .}$$

Obviously, it suffices to show that uniformly over all root sets $T$ in $G(n,p)$, $|T|=|R|$, the probability of such deviation is $o(\frac{1}{n^{|R|}})$. We use the same notation for $A_{i,j}$ as in the previous section. Due to Claim~\ref{cr:neigb} and the union bound, with probability $o(\frac{1}{n^{|R|}})$ for at least one of the constantly many sets $A_{i,j}$ in the decomposition of $T$ the number of common neigbours $\deg(A_{i, j})$ differs from $np^{k_i}$ by more than $\sqrt{2|R|np^{k_i}(1-p^ {k_i}) \ln n}$.  Let $\mathcal{S}_i$ be the set of all integers that differ from $np^{k_i}$ by at most $\sqrt{2|R|np^{k_i}(1-p^ {k_i}) \ln n}$. Then
$${\sf Pr}\left(|X(T) - Y(T)| > b_n \varepsilon\right) \leq \max\limits_{s_{i,j}\in\mathcal{S}_i} {\sf Pr}\left(|X(T) - Y(T)| > b_n \varepsilon \mid \deg(A_{i, j}) = s_{i,j}, \, i\in[r], j\in[m_i] \right) + o\left(\frac{1}{n^{|R|}}\right).$$
Let us first get an upper tail bound using the inequality from \cite[Proposition 2.44]{janson_ineq}. For convenience we recall this inequality below:

\begin{claim}[V. R\"{o}dl, A. Ruci\'{n}ski \cite{rodl_rucinski}]
\label{janson_up}
Let $\Gamma_p$ be a binomial random subset of a finite set $\Gamma$, and let $\mathcal{F}$ be a family of subsets in $\Gamma$. Let $Z=\sum_{F\in\mathcal{F}} I(F\subset\Gamma_p)$ count the number of times when $F\in\mathcal{F}$ appear as subsets of $\Gamma_p$. Let $D$ be the maximum (over $F$) number of sets in $\mathcal{F}$ that overlap with a single $F\in\mathcal{F}$. Then, for every $t\geq 0$,
  $${\sf Pr}(Z \geq {\sf E}Z + t) \leq (D + 1)\exp\left[{\frac{-t^2}{4(D + 1)({\sf E} Z + t/3)}}\right] .$$
\end{claim}

Now we fix $s_{i,j}\in\mathcal{S}_i$, $i\in[r], j\in[m_i]$, and also fix subsets $S_{i,j}\in[n]\setminus T$ of sizes $s_{i,j}$. Assume that  $N(A_{i,j}) = S_{i,j}$ for all $i\in[r], j\in[m_i]$. In order to apply Claim~\ref{janson_up}, we let $\Gamma$ to be the set of all edges that have both end-points outside $T$. Let $Z$ count the number of $(R,H)$-extensions of $T$. Then the family $\mathcal{F}$ consists of sets of edges induced by sets of vertices of size $O(n^{h-|R|})$, and thus $D=O(n^{h-|R|-2})$. Recall that $b_n = \Theta(n^{s+g-1} \sqrt{\frac{n}{\ln n}}) = \Theta(n^{h-|R|-1} \sqrt{\frac{n}{\ln n}})$ by (\ref{pattern_params_eq}). Therefore, from (\ref{both_bounds}) and the definition of $S$ it follows that that ${\sf E}(X(T) \left. \right| \deg(A_{1,1}), \ldots , \deg(A_{r, m_r})) = \Theta(n^{h-|R|})$. Thus, using Claim~\ref{janson_up}:
\begin{equation}\label{conc_up}
\begin{gathered}
   {\sf Pr}\left(X(T) - Y(T) > b(n)\varepsilon \left.  \right|  \deg(A_{i, j})=s_{i,j}, i\in[r], j\in[m_i] \right) \leq n^{O(1)}\exp\left[-\Theta\left(\frac{n}{\ln n}\right) \right] .
\end{gathered}
\end{equation}

To get the lower tail bound, we use the Janson's inequality \cite[Theorem 2.14]{janson_ineq}. Since the expected number $\overline{D}$ of edge-crossing extensions is $O(n^{2(h-|R|) - 2})$, we get:
\begin{equation}\label{conc_down}
    \begin{gathered}
    {\sf Pr}\left(X(T) - Y(T) < -b(n) \varepsilon \left. \right|  \deg(A_{i, j})=s_{i,j}, i\in[r], j\in[m_i] \right) \leq \exp\left[-\frac{b_n^2 \varepsilon^2}{2\overline{D}}\right] = \exp \left[-\Theta\left(\frac{n}{\ln n}\right)\right] .
\end{gathered}
\end{equation}

Combining (\ref{conc_up}) and (\ref{conc_down}), we finish the proof of (\ref{our_lemm}) and, thus, the proof of Theorem~\ref{th:main} as well.
\section{Further questions}
\label{sc:further}

We believe that our techniques can be used to prove the convergence of a rescaled maximum number of extensions even for non-symmetric $(R,H)$, while it should be hard to find the limit distribution.

\begin{figure}[htbp]
\centering
\includegraphics[scale=0.4]{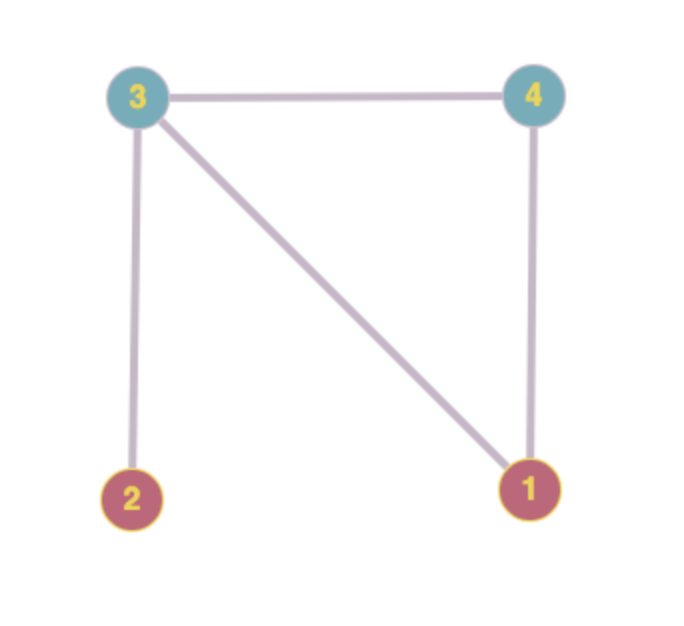}
\caption{a non-symmetric rooted graph, roots are in red}
\end{figure}
\label{fig:bad_example}

\noindent In particular, for the probably easiest non-symmetric $(R, H)$ consisting of two roots $v_1,v_2$ and two adjacent non-roots $u_1,u_2$ such that $u_1$ is adjacent to both $v_1,v_2$, and $u_2$ is only adjacent to $v_2$ (see Fig. \hyperref[fig:bad_example]{2}), we need a local limit theorem for vectors of dependent binomial random variables, which may be hard to eliminate.
Also, achieving a sufficient upper bound for $\Delta$ to apply Lemma~\ref{lm:formulation:good} could be technically very involved. Though we shall note that vertices of $H$ that are not adjacent to $R$ do not cause any additional difficulties.

Note that Bollob\'{a}s~\cite{bolobas}, Ivchenko~\cite{ivchenko} and Rodionov, Zhukovskii~\cite{common} studied also $m$-th maxima of cardinalities of common neighborhoods. It is of interest to get similar results for arbitrary symmetric extensions, while it might be not so evident when $r>1$ (let us recall that $r$ is the number of different cardinalities of root classes) or when $r=1$ and $m>2$.

Finally, our results can be generalised to $p=p(n)=o(1)$ (but $p>n^{-\varepsilon}$ for some small enough constant $\varepsilon>0$) when $r=1$. For larger $r$, the limit distribution that we get depends on $p$. So, for $r>1$ and $p=o(1)$, the limit behaviour  of the maximum number of extensions should be different.

\section*{Acknowledgements}

Stepan Vakhrushev is supported by Russian Science Foundation, project 22-11-00131.

\section*{Appendix}
\subsection*{A. Proof of Lemma~\ref{lm:our_tech}}

Let us denote $A_i=A_i(x):=\{Y_i > a_n+b_nx\}$, $B_i := \{X_i > a_n+b_nx\}$ for all $i \in [d]$. Note that it is sufficient to prove Lemma~\ref{lm:our_tech} for all $x\in\mathbb{R}$ such that $0<F(x)<1$. Let us fix such an $x\in\mathbb{R}$.  Find $\delta > 0$ such that $0 < F(x-\delta) \leq F(x+\delta) < 1$. Let $\varepsilon \in (0, \delta)$. We also denote $A_i^{\varepsilon} := A_i(x+\varepsilon)$.  The following inequalities hold:

$${\sf Pr}\left(\bigcup\limits_{i\in[d]} A_i^{\varepsilon}\right) - {\sf Pr}\left(\bigcup\limits_{i\in[d]} B_i\right) \leq {\sf Pr}\left(\bigcup\limits_{i\in[d]} A_i^{\varepsilon} \text{ }\backslash \bigcup\limits_{i\in[d]} B_i\right) \leq \sum\limits_{i\in[d]} {\sf Pr}(A_i^{\varepsilon} \backslash B_i) .$$

The condition (\ref{our_lemm}) implies $\sum\limits_{i\in[d]} {\sf Pr}(A_i^{\varepsilon} \backslash B_i) = o(1)$, so

 $${\sf Pr}\left(\bigcup\limits_{i\in[d]} B_i\right) \geq {\sf Pr}\left(\bigcup\limits_{i\in[d]} A_i^{\varepsilon}\right) - o(1) .$$
  
But from the first requirement in Lemma~\ref{lm:our_tech}

$$1 - {\sf Pr}\left(\bigcup\limits_{i\in[d]} A_i^{\varepsilon}\right) \xrightarrow{} F(x+\varepsilon) .$$

Recalling that $F$ is continuous and that the above holds for any $\varepsilon \in (0, \delta)$, we conclude that 

$$1-{\sf Pr}\left(\bigcup\limits_{i\in[d]} B_i\right) \leq F(x) + o(1) .$$

The lower bound $1-{\sf Pr}(\cup_{i\in[d]} B_i) \geq F(x) - o(1)$ is obtained similarly, using the events $A_i^{-\varepsilon}:=A_i(x-\varepsilon)$ and the relation $\sum\limits_{i\in[d]}{\sf Pr}(B_i \backslash A_i^{-\varepsilon}) = o(1)$ that follows directly from the condition (\ref{our_lemm}).

\subsection*{B. Proof of Lemma~\ref{lm:intersection}}
Since all the considered parameters are constants, it is sufficient to prove that, for any positive integers $k_1 \geq k_2 \geq k$ and $m_1, m_2$ whp the intersection of $U_{1, m_1}$ with $U_{2, m_2}$ does not equal to $k$. Let us denote this event by $A:=A(k_1, k_2, k, m_1, m_2)$. Let us separately consider the case when the second set is a subset of the first set, i.e. $k_1 > k_2 = k$.

\vspace{\baselineskip}

\begin{Proof1}

Let us estimate the probability of $A$ by the union bound over all choices of two sets $U_1\subset U_2$ on the role of $U_{1, m_1}$ and $U_{2, m_2}$:

\begin{equation}
\nonumber
\begin{gathered}
{\sf Pr}(A) \leq {n \choose k_1-k_2} \cdot {n-k_1+k_2 \choose k_2} \cdot {\sf Pr}(U_{1, m_1}=\{1, \ldots, k_1\}, U_{2, m_2}=\{1, \ldots, k_2\}).
\end{gathered}
\end{equation}

Fix $\varepsilon > 0$. From Theorem~\ref{th:max_com_neigb} the limit distribution of the maximum number of common neighbours implies that there exists a constant $C = C(\varepsilon)$ and an index $n_0$ starting from which:

\begin{equation}\label{system_of_restriction}
\nonumber
 \begin{cases}
   {\sf Pr}(|\Delta_{k_1, n}^{m_1} - a_{k_1,n}| \geq C \cdot \sqrt{\frac{n}{\ln n}}) < \varepsilon / 4  ,\\
  {\sf Pr}(|\Delta_{k_2, n}^{m_2} - a_{k_2,n}| \geq C \cdot  \sqrt{\frac{n}{\ln n}}) < \varepsilon / 4 .
 \end{cases}
\end{equation}

Hence, for $n > n_0$:

\begin{equation}
{\sf Pr}(A) \leq {n \choose k_1-k_2} \cdot {n-k_1+k_2 \choose k_2} \cdot {\sf Pr} \Bigg( |\deg(1, \ldots, k_i) - a_{k_i, n}| < C \sqrt{\frac{n}{\ln n}} \text{ for } i=1,2 \Bigg) + \frac{\varepsilon}{2} \hspace{3pt}.  \label{tech2}
\end{equation}

We write the internal probability in the following simple way: 

\begin{equation}
\begin{gathered}
{\sf Pr}\left( |\deg(1, \ldots, k_1) - a_{k_1, n}| < C \sqrt{\frac{n}{\ln n}}, |\deg(1, \ldots, k_2) - a_{k_2, n}| \leq C \sqrt{\frac{n}{\ln n}}\right) \leq \\ \leq \sum_{X \subset [n]: \frac{||X| - a_{k_1, n}|}{\sqrt{n/ \ln n}} \leq C} {\sf Pr}\Bigg( N(1, \ldots, k_1) = X \Bigg) {\sf Pr}\Bigg(\left.|\deg(1, \ldots, k_2) - a_{k_2, n} | \leq C \sqrt{\frac{n}{\ln n}}  \right|  N(1, \ldots, k_1) = X\Bigg) .
\end{gathered}
\label{int_prob}
\end{equation}

By the triangle inequality, the conditional probability in (\ref{int_prob}) is bounded from above by the probability that the number of neighbours of $U_2$ in $[n] \backslash (X \cup U_1)$ differs from $a_{k_2, n} - |X|$ by no more than $2C\sqrt{\frac{n}{\ln n}}$. By the de Moivre-Laplace limit theorem, the probability of the latter event approaches 0 as $n \to \infty$.

From (\ref{tech2}) and (\ref{int_prob}) we get:
\begin{equation}
\begin{gathered}
0 \leq {\sf Pr}(A)  \leq o\Big[ {n \choose k_1-k_2} \cdot {n-k_1+k_2 \choose k_2} \cdot {\sf Pr}\left(|\deg(1, \ldots, k_1) - a_{k_1, n}| \leq C \sqrt{\frac{n}{\ln n}}\right)\Big] + \frac{\varepsilon}{2} \hspace{3pt} .
\end{gathered}
\label{smthh}
\end{equation}

From (\ref{k_deg}) and (\ref{conv_exp}) it follows that ${\sf Pr}(|\deg(1, \ldots, k_1) - a_{k_1,n}| \leq C\sqrt{\frac{n}{\ln n}}) = O(n^{-k_1})$ implying that the first summand in the right hand side of (\ref{smthh}) approaches $0$ as $n \to \infty$. Due to arbitrariness of $\varepsilon$, the proof is completed.

\end{Proof1}

Now consider the case when none of the sets is nested in the other, i.e. $1 \leq k < \min(k_1, k_2)$. In \cite[Section 2.3.2]{common}, this statement is proven in the particular case $k_1 = k_2$. Our proof is similar, and we will use the bounds from  \cite[Section 2.3.2]{common} to get our results as well.
\vspace{\baselineskip}

\begin{Proof2}

First, let's narrow down the probability space to graphs with a ``small'' number of common neighbours:
$$
{\sf Pr}(A) \leq {\sf Pr}(A \cap \{G(n,p) \in \mathcal{Q}_n\}) + {\sf Pr}(G(n,p) \notin \mathcal{Q}_n).
$$

As discussed in Section 2, the second term tends to $0$. In what follows, we estimate only the joint probability. Fix $\varepsilon > 0$. From Theorem~\ref{th:max_com_neigb} there exists a constant $C = C(\varepsilon)$ such that starting from some $n_0 \in \mathbb{N}$:

\begin{equation}
\nonumber
 \begin{cases}
   {\sf Pr}(\Delta_{k_1, n}^{m_1} - a_{k_1,n} \leq -C \cdot \sqrt{\frac{n}{\ln n}}) < \varepsilon / 4, \\
  {\sf Pr}(\Delta_{k_2, n}^{m_2} - a_{k_2,n} \leq -C \cdot \sqrt{\frac{n}{\ln n}}) < \varepsilon / 4 .
 \end{cases}
\end{equation}
Then similarly to the previous case:
\begin{multline}
\nonumber
{\sf Pr}(A \cap \{G(n,p) \in \mathcal{Q}_n\}) \leq {n \choose k} {n - k \choose k_1 - k} {n - k_1 \choose k_2 - k} \times \\ \times {\sf Pr}\left(\frac{\deg(1, \ldots, k_1) - a_{k_1, n}}{\sqrt{n / \ln n}} > -C, \frac{\deg(k_1 - k + 1, \ldots, k_1 + k_2 - k) - a_{k_2, n}}{\sqrt{n / \ln n}} > -C, G(n,p) \in \mathcal{Q}_n \right) + \frac{\varepsilon}{2} .
\end{multline}

Hence, it suffices to prove that the fourth factor (probability of the event) is $o(n^{-(k_1 + k_2 - k)})$. Denote $b_1 = a_{k_1, n} - C \sqrt{\frac{n}{\ln n}}, b_2 = a_{k_2, n} - C \sqrt{\frac{n}{\ln n}}$.  It is obvious from the definition of $\mathcal{Q}_n$ that

\begin{multline}
{\sf Pr}\left(\deg([k_1]) > b_1, \deg([k_1+k_2]\setminus [k_1] - k) > b_2, G(n,p) \in \mathcal{Q}_n \right) \leq \\ \leq \sum\limits_i {\sf Pr}(\xi_{n, p^k}=i) {\sf Pr}(\xi_{i, p^{k_1-k}} > b_1 - (k_2-k)) {\sf Pr}(\xi_{i, p^{k_2-k}} > b_2 - (k_1-k)) + \\ + {\sf Pr}(\xi_{n, p^k} \leq np^k-\sqrt{2(k_1+k_2)p^k(1-p^k)n \ln n }),
\label{fourth_factor}
\end{multline}
\noindent
where the summation is over $i\in \left(np^k-\sqrt{2(k_1+k_2)p^k(1-p^k)n \ln n}, \Gamma_k\right]$.  From Claim~\ref{cr:neigb} we get that the second term is $\frac{n^{-(k_1 + k_2)}(1 + o(1))}{2\sqrt{(k_1 + k_2)\pi \ln n}} = o(n^{-(k_1+k_2)})$. Therefore, it suffices to estimate only the first sum.

By the de Moivre--Laplace limit theorem, uniformly over $i$:
\begin{equation}\label{first_var}
    {\sf Pr}(\xi_{n, p^k} = i) = \frac{\exp \left[{-\frac{(np^k - i)^2}{2np^k(1-p)^k}}\right]}{\sqrt{2\pi n p^k(1-p^k)}}(1+ o(1)) .
\end{equation}

By the de Moivre--Laplace limit theorem (here we skip the computations, that can be found in~\cite[Section 2.3.2]{common}):

\begin{equation}{\label{second_var}}
    {\sf Pr}(\xi_{i, p^{k_1 - k}} > b_1 - (k_2 - k)) \leq \frac{\sqrt{1 - p^{k_1-k}}e^{-\frac{(b_1 - i p^{k_1-k})^2}{2i p^{k_1-k}(1-p^{k_1-k})}}(1+o(1))}{\sqrt{2\pi \ln n}\left(\sqrt{2k_1(1-p^{k_1})} - \sqrt{2k(p^{k_1-k} - p^{k_1})} \right)} ,
\end{equation}
\noindent
and the same bound holds true with $k_1$ replaced with $k_2$ and $b_1$ replaced with $b_2$. From (\ref{first_var}) and (\ref{second_var}), we get that the first summand in right-hand side of (\ref{fourth_factor}) is $\frac{O(1)}{\sqrt{n} \ln n} \sum e^{-{g(i)}}$
, where

$$g(i) = \frac{(np^k - i)^2}{2np^k(1-p^k)} + \frac{(ip^{k_1 - k} - b_1)^2}{2ip^{k_1 - k}(1-p^{k_1 - k})} + \frac{(ip^{k_2 - k} - b_2)^2}{2ip^{k_2 - k}(1-p^{k_2 - k})} .$$

Denote $i = np^k + x\sqrt{np^k(1-p^k)\ln n}$, $x \in (-\sqrt{2(k_1 + k_2)}, \sqrt{2k}]$. Then the first term in $g(i)$ becomes $\frac{x^2}{2} \ln n$.
After the replacement, we get:
\begin{equation}
\nonumber
    g(i) = \tilde{g}_p(x) \ln n + \hat{g}_p(x) \ln \ln n(1 + o(1)),
\end{equation}
where 

\begin{equation}
\nonumber
\begin{gathered}
\tilde{g}_p(x) = \frac{x^2}{2} + \frac{x^2(p^{k_1 - k} - p^{k_1})}{2(1 - p^{k_1 - k})} + \frac{x^2(p^{k_2 - k} - p^{k_2})}{2(1 - p^{k_2 - k})} - \frac{ 2\sqrt{2k_1}\sqrt{(p^{k_1 - k} - p^{k_1})(1 - p^{k_1})}x}{2(1-p^{k_1 - k})} - \\ - \frac{ 2\sqrt{2k_2}\sqrt{(p^{k_2 - k} - p^{k_2})(1 - p^{k_2})}x}{2(1-p^{k_2 - k})} + \frac{2k_1(1-p^{k_1})}{2(1-p^{k_1 - k})} + \frac{2k_2(1-p^{k_2})}{2(1-p^{k_2 - k})},
\end{gathered}
\end{equation}

\noindent
and $\hat{g}_p(x)$ is negative and bounded from below by a constant (in the same way as in~\cite[Section 2.3.2]{common}). It follows from the size of the summation segment that it suffices for us to show that $\tilde{g}_p(x) \geq k_1+k_2-k + \omega(\frac{\ln \ln n}{\ln n})$. We need the positive term $\omega(\frac{\ln \ln n}{\ln n})$ to overcome the negative contribution of $\hat g_p(x)$.

We set $\tilde g_p(x)=\frac{1}{2}(\tilde g_{1,p}(x)+\tilde g_{2,p}(x))$, where
$$
 \tilde g_{j,p}(x)=\frac{x^2(1+p^{k_j-k} - 2p^{k_j} )- 4\sqrt{2k_j}\sqrt{(p^{k_j-k} - p^{k_j})(1-p^{k_j})}x + 4k_j(1-p^{k_j})}{2(1-p^{k_j-k})}, \quad j=1,2.
$$

In the same way as in~\cite[Section 2.3.2]{common}, we get that, for every $j\in\{1, 2\}$, $\tilde{g}_{j,p}(x) \geq 2k_j - k + \omega(\frac{\ln \ln n}{\ln n})$ completing the proof.

\end{Proof2}
\end{document}